\documentclass[11pt]{amsart}
\usepackage{amsmath}
\usepackage{amssymb}

\begin{document}


\newtheorem{theorem}{Theorem}
\newtheorem{example}[theorem]{Example}
\newtheorem{definition}[theorem]{Definition}

\newcommand{\Eqref}[1]{(\ref{#1})}



\title{Transmutation operators: construction and applications}



\author[V. V. Kravchenko]{Vladislav V. Kravchenko}

\address[V. V. Kravchenko and S. M. Torba]{Department of Mathematics, CINVESTAV del IPN, Unidad Querétaro,\\
Libramiento Norponiente No. 2000, Fracc. Real de Juriquilla\\
Querétaro, Qro. C.P. 76230 MEXICO}

\email[V. V. Kravchenko]{vkravchenko@math.cinvestav.edu.mx}

\author[S. M. Torba]{Sergii M. Torba}

\email[S. M. Torba]{storba@math.cinvestav.edu.mx}

\author[K. V. Khmelnytskaya]{Kira V. Khmelnytskaya}

\address[K. V. Khmelnytskaya]{Faculty of Engineering, Autonomous University of
Quer\'{e}taro,  M\'{e}xico}

\email[K. V. Khmelnytskaya]{khmel@uaq.edu.mx}

\thanks{This work has been supported by CONACYT, Mexico via the projects 166141 and 222478.}

\dedicatory{ }

\subjclass[2010]{34A25, 34B24, 35C05, 35J05, 34L40, 42C10, 65L05, 65L70, 65M80}

\keywords{Sturm-Liouville equation, Transmutation operator, Neumann series of Bessel functions, Spectral problem, Complete family of solutions, Method of fundamental solutions}

\begin{abstract}
Recent results on the construction and applications of the
transmutation (transformation) operators are discussed. Three new
representations for solutions of the one-dimensional Schr\"{o}dinger equation
are considered. Due to the fact that they are obtained with the aid of the
transmutation operator all the representations possess an important for
practice feature. The accuracy of the approximate solution is independent of
the real part of the spectral parameter. This makes the representations
especially useful in problems requiring computation of large sets of eigendata
with a nondeteriorating accuracy.

Applications of the exact representations for the transmutation
operators to partial differential equations are discussed as well. In
particular, it is shown how the methods based on complete families of
solutions can be extended onto equations with variable coefficients.
\end{abstract}
\maketitle

\section{Transmutation operators}

Transmutation operators also called transformation operators are a widely used
tool in the theory of linear differential equations (see, e.g., \cite{Gilbert}%
, \cite{Carroll}, \cite{Colton} \cite{LevitanInverse}, \cite{Marchenko},
\cite{Sitnik} and many other publications). In particular, let $q\in C[-b,b]$
be a complex valued function. Consider the Sturm-Liouville equation
\begin{equation}
Ay:=y^{\prime\prime}-q(x)y=-\omega^{2}y.\label{S-L}%
\end{equation}
It is well known (see, e.g., \cite{Marchenko}) that there exists a Volterra
integral operator $T$ called the transmutation (or transformation) operator
defined on $C[-b,b]$ by the formula
\[
Tu(x)=u(x)+\int_{-x}^{x}K(x,t)u(t)dt
\]
such that for any $u\in C^{2}[-b,b]$ the following equality is valid%
\[
ATu=Tu^{\prime\prime}%
\]
and hence any solution of (\ref{S-L}) can be written as $y=T[u]$ where
$u(x)=c_{1}\cos\omega x+c_{2}\sin\omega x$ with $c_{1}$ and $c_{2}$ being
arbitrary constants.

The transmutation kernel $K$ is a solution of a certain Goursat problem for
the hyperbolic equation
\[
\left(  \frac{\partial^{2}}{\partial x^{2}}-q(x)\right)  K(x,t)=\frac
{\partial^{2}}{\partial t^{2}}K(x,t).
\]

\section{Construction of the transmutation kernel and new representations for
solutions of the Sturm-Liouville equation}

In spite of fundamental importance of the transmutation kernel $K$ in the
theory of linear differential equations, besides the method of successive
approximations derived directly from the Goursat problem (see, e.g.,
\cite{Colton}) very few attempts of its practical construction have been
reported. In this relation we mention the paper \cite{Boumenir2006} where
analytic approximation formulas for the integral kernel were obtained and the
recent publications \cite{KT Transmut}, \cite{KT AnalyticApprox} where another
procedure of analytical approximation was proposed.

To the difference of those previous results, in the recent paper \cite{KNT} an
exact representation for $K$ in the form of a Fourier-Legendre series with
explicit formulas for the coefficients was obtained. Suppose that $q\in
W_{2}^{-1}[-b,b]$ (that is, $q$ can be a piecewise continuous function, may
have a singularity, e.g., $q(x)\sim c/x$, etc.). In this case $K(x,t)$ is an
$L_{2}$-function. Under these conditions the following theorem was proved in
\cite{KNT}.

\begin{theorem}
The kernel $K$ has the form
\begin{equation}
K(x,t)=\sum_{n=0}^{\infty}\frac{\beta_{n}(x)}{x}P_{n}\left(  \frac{t}%
{x}\right)  \label{K}%
\end{equation}
where for every $x\in\lbrack-b,b]$ the series converges with respect to $t$ in
the $L_{2}$-norm (if $q\in C[-b,b]$ the series converges uniformly),
\[
\beta_{n}(x)=\frac{2n+1}{2}\biggl(\sum_{k=0}^{n}\frac{l_{k,n}\varphi_{k}%
(x)}{x^{k}}-1\biggr),
\]
with $l_{k,n}$ being the coefficient at $x^{k}$ of the Legendre polynomial
$P_{n}$, and $\varphi_{k}$ being the so-called formal powers\ constructed as
follows (see \cite{KrCV08}, \cite{KrPorter2010}).
\end{theorem}

\begin{definition}[Formal powers $\varphi_{k}$] Let $f$ be a solution of
\begin{equation}
f^{\prime\prime}-q(x)f=0,\quad x\in\lbrack-b,b],\label{SLhom}%
\end{equation}%
\[
f(0)=1,\quad f^{\prime}(0)=0.
\]
Then $\left\{  \varphi_{k}\right\}  _{k=0}^{\infty}$ are defined by the
equalities%
\[
\varphi_{k}=
\begin{cases}
fX^{(k)}, & k \text{ odd,}\\
f\widetilde{X}^{(k)}, & k \text{ even},
\end{cases}
\]
where%
\[
X^{(0)}\equiv1,\qquad X^{(n)}(x)=n\int_{0}^{x}X^{(n-1)}(s)\left(
f^{2}(s)\right)  ^{(-1)^{n}}\,\mathrm{d}s,
\]
and%
\[
\widetilde{X}^{(0)}\equiv1,\qquad\widetilde{X}^{(n)}(x)=n\int_{0}%
^{x}\widetilde{X}^{(n-1)}(s)\left(  f^{2}(s)\right)  ^{(-1)^{n-1}}%
\,\mathrm{d}s.
\]

\end{definition}

It is worth mentioning that $\varphi_{k}$ are easily computable (at least numerically) in practice (see, e.g., \cite{KNT} for additional details).

A representation for the kernel K leads to a representation for the solution
of (\ref{S-L}). Let $u(\omega,x)$ denote the solution of (\ref{S-L})
satisfying the initial conditions%
\[
u(\omega,0)=1,\qquad u^{\prime}(\omega,0)=i\omega.
\]
Then we have %

\[
u(\omega,x)=e^{i\omega x}+\int_{-x}^{x}K(x,t)e^{i\omega t}dt.
\]
Substitution of (\ref{K}) into the last integral gives us the equality
\cite{KNT}
\[
u(\omega,x)=e^{i\omega x}+\sum_{n=0}^{\infty}\beta_{n}(x)\int_{-1}^{1}%
P_{n}\left(  y\right)  e^{i\omega xy}dy=e^{i\omega x}+\sum_{n=0}^{\infty}%
i^{n}\beta_{n}(x)j_{n}(\omega x)
\]
where $j_{n}(z)=\sqrt{\frac{\pi}{2z}}J_{n+1/2}(z)$ are spherical Bessel
functions. The series converges uniformly with respect to $x$.\textbf{\ }

Moreover, take $\omega\in\mathbb{R}$. Consider
\[
K_{N}(x,t)=\sum_{n=0}^{N}\frac{\beta_{n}(x)}{x}P_{n}\left(  \frac{t}%
{x}\right)
\]
and
\[
u_{N}(\omega,x)=e^{i\omega x}+\sum_{n=0}^{N}i^{n}\beta_{n}(x)j_{n}(\omega
x),\quad x>0.
\]
We have \cite{KNT}
\begin{align*}
\left\vert u(\omega,x)-u_{N}(\omega,x)\right\vert  & =\left\vert \int_{-x}%
^{x}\left(  K(x,t)-K_{N}(x,t)\right)  e^{i\omega t}dt\right\vert \\
& \leq\left\Vert K(x,\cdot)-K_{N}(x,\cdot)\right\Vert _{L_{2}(-x,x)}%
\,\left\Vert e^{i\omega t}\right\Vert _{L_{2}(-x,x)}\\
& =\varepsilon_{N}(x)\sqrt{2x}
\end{align*}
---independent of $\omega$.
More generally, for any $\omega\in\mathbb{C}$, $\omega\neq0$ belonging to the
strip $\left\vert \operatorname{Im}\omega\right\vert \leq C$, $C\geq0$,
\[
\left\vert u(\omega,x)-u_{N}(\omega,x)\right\vert \leq\varepsilon_{N}%
(x)\frac{\sinh(Cx)}{C}.
\]
This $\omega$-independence of the approximation accuracy was shown in
\cite{KNT} to give a very fast and efficient method for computing large sets
of eigendata with a nondeteriorating accuracy. In \cite{Kr2} it was
generalized onto perturbed Bessel equations, and in \cite{KT2016} onto
Sturm-Liouville equations.

Another representation for the kernel $K$ and as a corollary for the solutions
of (\ref{S-L}) was obtained in \cite{Kr3}. Consider the following extension of
the transmutation kernel $K$,
\[
\widetilde{K}(x,t)=
\begin{cases}
K(x,t), & -x\leq t\leq x,\\
0, & -\infty<t<-x.%
\end{cases}
\]
Then
\[
u(\omega,x)=e^{i\omega x}+\int_{-\infty}^{x}\widetilde{K}(x,y)e^{i\omega
y}dy=e^{i\omega x}\left(  1+\int_{0}^{\infty}\widetilde{K}(x,x-t)e^{-i\omega
t}dt\right)  .
\]
Consider
\[
\widetilde{K}(x,x-t)=\mathbf{k}(x,t)e^{-t}.
\]
The function $\mathbf{k}(x,\cdot)$ then belongs to the space $L_{2}\left(
0,\infty;e^{-t}\right)  $ equipped with the scalar product $\left\langle
u,v\right\rangle :=\int_{0}^{\infty}u(t)\overline{v}(t)e^{-t}dt$. Thus,
$\mathbf{k}(x,\cdot)$ admits a Fourier-Laguerre expansion convergent in the
corresponding norm,
\[
\mathbf{k}(x,t)=\sum_{n=0}^{\infty}a_{n}(x)L_{n}(t).
\]
The kernel has the form \cite{Kr3}%
\[
\widetilde{K}(x,y)=\sum_{n=0}^{\infty}a_{n}(x)L_{n}(x-y)e^{-(x-y)},
\]
with the coefficients $a_{n}$ defined by%
\begin{equation}
a_{n}(x)=\sum_{j=0}^{n}\left(  -1\right)  ^{j}\left(  \varphi_{j}%
(x)-x^{j}\right)  \sum_{k=j}^{n}\left(  -1\right)  ^{k}\frac{n!}%
{(n-k)!k!(k-j)!j!}x^{k-j}.\label{a_n}%
\end{equation}
The solution $u(\omega,x)$ has the form \cite{Kr3}%
\begin{equation}
u(\omega,x)=e^{i\omega x}\left(  1+\sum_{n=0}^{\infty}a_{n}(x)\frac{\left(
i\omega\right)  ^{n}}{\left(  1+i\omega\right)  ^{n+1}}\right)  .\label{u}%
\end{equation}

The following estimate is valid for any $\omega\in\mathbb{R}$,
\begin{equation}
\left\vert u(\omega,x)-u_{N}(\omega,x)\right\vert \leq\varepsilon_{N}\left(
x\right)  ,\label{estimate omega}%
\end{equation}
where
\[
u_{N}(\omega,x):=e^{i\omega x}\left(  1+\sum_{n=0}^{N}a_{n}(x)\frac{\left(
i\omega\right)  ^{n}}{\left(  1+i\omega\right)  ^{n+1}}\right)  ,
\]
and $\varepsilon_{N}\left(  x\right)  $ is a nonnegative function independent
of $\omega$ and such that $\varepsilon_{N}\left(  x\right)  \rightarrow0$ for
all $x\in\left[  -b,b\right]  $ when $N\rightarrow\infty$. More generally,%
\[
\left\vert u(\omega,x)-u_{N}(\omega,x)\right\vert \leq\frac{\varepsilon
_{N}(  x)  e^{-\operatorname{Im}\omega x}}{\sqrt
{1-2\operatorname{Im}\omega}},\qquad\text{when }\operatorname{Im}\omega<1/2.
\]

Consideration of another extension of the transmutation kernel defined by
\[
\widetilde{K}(x,y):=
\begin{cases}
K(x,y) & \text{when  }x\in\left[  -b,b\right]  \text{ and } y\in\left[  -x,x\right]  \\
0 & \text{otherwise}.
\end{cases}
\]
leads to the following series expansion
\[
\widetilde{K}(x,y)=\sum_{n=0}^{\infty}c_{n}(x)H_{n}(y)e^{-y^{2}}%
\]
where $H_{n}$ stands for an Hermite polynomial of order $n$ and the
coefficients $c_{n}$ are to be found. Note that
\[
\int_{-\infty}^{\infty}\widetilde{K}(x,y)H_{n}(y)dy=\sqrt{\pi}n!2^{n}c_{n}(x).
\]
Hence
\begin{equation}
c_{n}(x)=\frac{1}{\sqrt{\pi}n!2^{n}}\int_{-x}^{x}K(x,y)H_{n}(y)dy=\frac
{1}{\sqrt{\pi}n!2^{n}}\sum_{k=0}^{n}h_{k,n}\left(  \varphi_{k}(x)-x^{k}%
\right)  \label{cn}%
\end{equation}
where $h_{k,n}$ denotes the coefficient of $x^{k}$ from the Hermite polynomial
$H_{n}(x)$.

This leads to another representation for the solution of (\ref{S-L}),
\begin{align*}
u(\omega,x) &  =e^{i\omega x}+\int_{-\infty}^{\infty}\widetilde{K}%
(x,y)e^{i\omega y}dy\\
&  =e^{i\omega x}+\sum_{n=0}^{\infty}c_{n}(x)\int_{-\infty}^{\infty}%
H_{n}(y)e^{i\omega y}e^{-y^{2}}dy\\
&  =e^{i\omega x}+\sqrt{\pi}e^{-\frac{\omega^{2}}{4}}\sum_{n=0}^{\infty}%
c_{n}(x)\left(  i\omega\right)  ^{n}.
\end{align*}
Consider the partial sum
\[
u_{N}(\omega,x)=e^{i\omega x}+\sqrt{\pi}e^{-\frac{\omega^{2}}{4}}\sum
_{n=0}^{N}c_{n}(x)\left(  i\omega\right)  ^{n}.
\]
Then it is easy to see that%
\[
\left\vert u(\omega,x)-u_{N}(\omega,x)\right\vert \leq\pi^{\frac{1}{4}%
}e^{\frac{\left(  \operatorname{Im}\omega\right)  ^{2}}{2}}\varepsilon_{N}(x)
\]
which means that the truncation error is uniformly bounded in any strip
$\left\vert \operatorname{Im}\omega\right\vert \leq C$.

\section{Applications to PDEs}

Exact representations of the transmutation kernel lead to numerous
applications for partial differential equations admitting certain symmetry. In particular, let us consider the possibility to obtain complete systems of
solutions. For example, application of the transmutation operator $T$ to a
complete system of harmonic functions leads to a complete system of solutions
of the equation
\begin{equation}
\left(  \Delta-q(x)\right)  u(x,y)=0.\label{Helm}%
\end{equation}
Indeed,
\[
\left(  \Delta-q(x)\right)  T=T\Delta
\]
whenever the domain of interest is such that the integration in $T$ is well defined.

\textbf{Example} Harmonic polynomials ($\operatorname{Re}z^{n}$ and
$\operatorname{Im}z^{n}$) can be written in the form
\begin{align*}
p_{0}(x,y) &  =1,\\
p_{2m+1}(x,y) &  =\operatorname{Re}z^{m+1}=\sum_{\text{even }k=0}^{m}
\left(  -1\right)  ^{\frac{k}{2}}\binom{m+1}{k}x^{m+1-k}y^{k},\qquad m\ge 0,\\
p_{2m}(x,y) &  =\operatorname{Re}\left(  iz^{m}\right)  =\sum_{\text{odd }k=1}^{m}\left(  -1\right)  ^{\frac{k+1}{2}}\binom{m}{k}x^{m-k}y^{k},\qquad m\ge 1.
\end{align*}

Since
\[
T:x^{k}\mapsto\varphi_{k}(x),
\]
the following functions are the images of $p_{m}$ under the action of $T$ and
represent a complete system of solutions of (\ref{Helm})
\begin{align*}
u_{0}(x,y) &  =f(x),\\
u_{2m+1}(x,y) &  =\sum_{\text{even }k=0}^{m+1}
\left(  -1\right)  ^{\frac{k}{2}}\binom{m+1}{k}\varphi_{m+1-k}(x)y^{k},\qquad m\ge 0,\\
u_{2m}(x,y) &  =\sum_{\text{odd }k=1}^{m}
\left(  -1\right)  ^{\frac{k+1}{2}}\binom{m}{k}\varphi_{m-k}(x)y^{k},\qquad m\ge 1.
\end{align*}
Other complete systems of solutions can be obtained.

\textbf{Example }The method of fundamental solutions (discrete sources) (see,
e.g., \cite{Alexidze}, \cite{Doicu et al}) can be extended onto equations with
variable coefficients. Consider the fundamental solution for the Laplace
operator on the plane%

\[
\log\left\vert x+iy-(\eta+i\xi)\right\vert =\log\left\vert x-Z\right\vert .
\]
Application of $T$ leads to the following integrals%
\[
T\left[  \log\left\vert x-Z\right\vert \right]  =\log\left\vert x-Z\right\vert
+\sum_{n=0}^{\infty}\frac{\beta_{n}(x)}{x}\int_{-x}^{x}P_{n}\left(  \frac
{t}{x}\right)  \log\left\vert t-Z\right\vert \,dt.
\]
Their calculation gives us the image of the fundamental solution in the form
\begin{align*}
 T\left[  \log\left\vert x-Z\right\vert \right]  &=\log\left\vert
x-Z\right\vert +\beta_{0}(x)\operatorname{Re}\left(  \log\left(  \left(
Z+x\right)  \left(  Z-x\right)  \right)  +2Q_{1}\left(  \frac{Z}{x}\right)
\right)  \\
&\quad +2\sum_{n=1}^{\infty}\frac{\beta_{n}(x)}{2n+1}\operatorname{Re}\left(
Q_{n+1}\left(  \frac{Z}{x}\right)  -Q_{n-1}\left(  \frac{Z}{x}\right)
\right)
\end{align*}
where $Q_{n}$ are Legendre functions of the second kind.

Similar considerations can be applied to systems arising in hypercomplex
analysis (see, e.g., \cite{GS2}) and may lead to extensions of well known
methods based on monogenic polynomials or other complete systems of solutions
onto systems with variable coefficients \cite{AQA}, \cite{APFT} important,
e.g., in electromagnetic theory and quantum physics.


\begin{thebibliography}{99}

\bibitem {Alexidze}{\sc M.\,A.\,Alexidze}, \emph{Fundamental Functions in Approximate
Solutions of Boundary Value Problems} (in Russian), Nauka, Moscow, 1991.

\bibitem {Gilbert}{\sc H.\,Begehr and R.\,Gilbert}, \emph{Transformations,
transmutations and kernel functions, vol. 1--2}, Longman Scientific \&
Technical, Harlow, 1992.

\bibitem {Boumenir2006}{\sc A.\,Boumenir}, \emph{The approximation of the
transmutation kernel}, {\rm J. Math. Phys.} {\bf 47} (2006), 013505.

\bibitem {Carroll}{\sc R.\,W.\,Carroll}, \emph{Transmutation theory and applications, Mathematics Studies, Vol. 117}, North-Holland, 1985.

\bibitem {Colton}{\sc D.\,Colton}, \emph{Solution of boundary value problems by the method of integral operators}, Pitman, London, 1976.

\bibitem {Doicu et al}{\sc A.\,Doicu, Yu.\,Eremin and Th.\,Wriedt}, \emph{Acoustic and
Electromagnetic Scattering Analysis,} Acad. Press, London, 2000.

\bibitem {GS2}{\sc K.\,G\"{u}rlebeck and W.\,Spr\"{o}ssig}, \emph{Quaternionic and
Clifford Calculus for Physicists and Engineers}, John Wiley \&
Sons, Chichester, 1997.

\bibitem {KrCV08}{\sc V.\,V.\,Kravchenko}, \emph{A representation for solutions of the
Sturm-Liouville equation}, {\rm Complex Var. Elliptic Equ.} {\bf 53} (2008) 775--789.

\bibitem {Kr3}{\sc V.\,V.\,Kravchenko}, \emph{Construction of a transmutation for the
one-dimensional Schr\"{o}dinger operator and a representation for solutions}, submitted, available from arXiv: 1612.09577.

\bibitem {AQA}{\sc V.\,V.\,Kravchenko}, \emph{Applied quaternionic analysis},
Heldermann Verlag, Lemgo, 2003.

\bibitem {APFT}{\sc V.\,V.\,Kravchenko}, \emph{Applied pseudoanalytic function theory, Series: Frontiers in Mathematics},
Birkh\"{a}user, Basel, 2009.

\bibitem {KrPorter2010}{\sc V.\,V.\,Kravchenko and R.\,M.\,Porter}, \emph{Spectral
parameter power series for Sturm-Liouville problems}, {\rm Math. Methods Appl. Sci.}
{\bf 33} (2010) 459--468.

\bibitem {KNT}{\sc V.\,V.\,Kravchenko, L.\,J.\,Navarro and S.\,M.\,Torba},
\emph{Representation of solutions to the one-dimensional Schr\"{o}dinger
equation in terms of Neumann series of Bessel functions}, {\rm Appl. Math. Comput.} {\bf 314} (2017) 173--192.

\bibitem {KT Transmut}{\sc V.\,V.\,Kravchenko and S.\,M.\,Torba}, \emph{Construction of
transmutation operators and hyperbolic pseudoanalytic functions}, {\rm Complex
Anal. Oper. Theory} {\bf 9} (2015) 389--429.

\bibitem {KT AnalyticApprox}{\sc V.\,V.\,Kravchenko and S.\,M.\,Torba}, \emph{Analytic
approximation of transmutation operators and applications to highly accurate
solution of spectral problems}, {\rm J. Comput. Appl. Math.} {\bf 275} (2015) 1--26.

\bibitem {KT2016}{\sc V.\,V.\,Kravchenko and S.\,M.\,Torba}, \emph{A Neumann series of
Bessel functions representation for solutions of Sturm-Liouville equations},
submitted, available from: arXiv:1612.08803.

\bibitem {Kr2}{\sc V.\,V.\,Kravchenko, S.\,M.\,Torba and R.\,Castillo-Perez}, \emph{ A Neumann
series of Bessel functions representation for solutions of perturbed Bessel
equations}, {\rm Appl. Analysis} {\bf } (2017), 28 pp. Published online, doi:10.1080/00036811.2017.1284313.

\bibitem {LevitanInverse}{\sc B.\,M.\,Levitan}, \emph{Inverse Sturm-Liouville
problems}, VSP, Zeist, 1987.

\bibitem {Marchenko}{\sc V.\,A.\,Marchenko}, \emph{Sturm-Liouville operators and
applications}, Birkh\"{a}user, Basel, 1986.

\bibitem {Sitnik}{\sc S.\,M.\,Sitnik}, \emph{Transmutations and applications: a
survey}, arXiv:1012.3741v1, originally published in the book: \emph{Advances
in Modern Analysis and Mathematical Modeling}, Editors: {\sc Yu.\,F.\,Korobeinik and
A.\,G.\,Kusraev}, Vladikavkaz Scientific Center of the Russian Academy
of Sciences and Republic of North Ossetia--Alania, Vladikavkaz, 2008, 226--293.
\end{thebibliography}
\end{document}